\documentclass{iopart}
\usepackage{iopams}
\usepackage{geometry}                % See geometry.pdf to learn the layout options. There are lots.
\usepackage{graphicx}
\usepackage{epstopdf}
\usepackage{amsthm}

\usepackage{cite}

\usepackage[pagebackref,pdfusetitle]{hyperref}

\usepackage{url}

\usepackage[percent]{overpic}

\usepackage{color}

\newcommand{\morespacearray}{\renewcommand{\arraystretch}{1.5}}

\newtheorem{proposition}{Proposition}

\newtheorem{theorem}[proposition]{Theorem}

\begin{document}

%\title{On the preservation of certain quadratic integrals \\ \hspace{3cm}using the Kahan discretization}

\title{Geometric and integrability properties of Kahan's method:\\ \hspace{1.5cm}The preservation of certain quadratic integrals}

\author{E Celledoni$^1$, D I McLaren$^2$, B Owren$^1$ and G R W Quispel$^2$}

\address{$^1$ 	Department of Mathematical Sciences,
	NTNU,
	7491 Trondheim,
	Norway\ead{elenac@ntnu.no}\ead{brynjulf.owren@ntnu.no}}
%\address{$^2$ 	Institute of Fundamental Sciences,
%	Massey University,
%	Private Bag 11 222, Palmerston North 4442, New Zealand\ead{r.mclachlan@massey.ac.nz}}
\address{$^2$ 	Department of Mathematics,
	La Trobe University,
	Bundoora, VIC 3083, Australia   \ead{d.mcLaren@latrobe.edu.au}\ead{r.quispel@latrobe.edu.au}}

\begin{abstract}
\noindent
Given a quadratic vector field on $\mathbb{R}^n$ possessing a quadratic first integral depending on two of the independent variables, we give a constructive proof that Kahan's discretization method exactly preserves a nearby modified integral. Building on this result, we present a  family of integrable quadratic vector fields (including the Euler top) whose Kahan discretization is a family of integrable maps. 
\end{abstract}

\section{Introduction}
Most ordinary differential equations (ODEs) cannot be solved exactly in closed form. In general, the next best thing is to use a numerical integration method that preserves one or more geometric properties (or a nearby property) of a given ODE exactly. Indeed, the field of geometric numerical integration is devoted to this endeavour. Thus, methods have been developed that preserve either symplectic structure, or energy and other integrals, or phase space volume, or (reversing) symmetries, or dissipation, etc exactly.
However, when it comes to preserving two or more of the properties, in many cases it may be unknown or impossible to do this \cite{HLW}. A notable exception to this situation is presented by Kahan's method, which in many cases preserves both (modified versions of) all integrals as well as the volume form exactly. Kahan introduced his ``unconventional" discretization method of quadratic vector fields in 1993 \cite{K}. Given a quadratic vector field in $\mathbb{R}^n$:
\begin{equation} \label{quadvectorfield}
\frac{\mathrm{d}x_i}{\mathrm{d}t} = \sum_{j,k} a_{ijk} x_j x_k + \sum_j b_{ij} x_j + c_i,\quad i=1,\ldots,n.
\end{equation}
Kahan's discretization is given by
\begin{equation} \label{Kahansmethod}
   \frac{x_i'-x_i}{h} = \sum_{j,k} a_{ijk} \frac{x_j'x_k+x_jx_k'}{2} + \sum_j  b_{ij}\frac{x_j+x_j'}{2} + c_i,\quad i=1,\ldots,n
\end{equation}
where $h$ denotes the discrete time step and
$$
      x_i \approx x_i(mh);\ x_i' \approx x_i((m+1)h).
$$
It is important to note that Kahan's method (\ref{Kahansmethod}), as well as its inverse, are linearly implicit, and hence define a birational map.

Kahan's method was independently rediscovered by Hirota and Kimura \cite{HK,KH}, and in 2011 Petrera, Pfadler and Suris \cite{PPS} applied  Kahan's method to a large number of integrable quadratic vector fields and showed that the discretization in most cases preserved the integrability.
In \cite{CMOQ} and \cite{CMMOQ}, we have shown that all linear and cubic integrals preserved in \cite{PPS} using Kahan's method can be explained by general propositions unrelated to complete integrability. 

In the present paper, Theorem~\ref{theo-1} similarly yields the preservation of many (but not all) quadratic integrals, and building on this, Theorem~\ref{theo-2} gives a 10-parameter family of integrable maps in $\mathbb{R}^3$.

\section{The preservation of quadratic integrals of two variables by Kahan's method}
Let the quadratic ODE
\begin{equation} \label{GeneralODE}
   \frac{\mathrm{d}\mathbf{x}}{\mathrm{d}t} = f(\mathbf{x}),\qquad \mathbf{x}\in\mathbb{R}^n
\end{equation}
possess a quadratic integral $I$ in 2 variables. W.l.o.g. we can choose the variables to be $x_1$ and $x_2$ with
\begin{equation} \label{Integral2vars}
I(x_1,x_2)= \frac12 a_1x_1^2 + a_2 x_1 x_2 + \frac12 a_3 x_2^2 + a_4 x_1 + a_5 x_2.
\end{equation}
It follows that the first two components of the vector field can be written
%\begin{equation} \label{Firsttwocomp}
%\begin{split}
%\frac{\mathrm{d}x_1}{\mathrm{d}t} & = A(x)\,\frac{\partial I}{\partial x_2} \\
%\frac{\mathrm{d}x_2}{\mathrm{d}t} & = -A(x)\,\frac{\partial I}{\partial x_1} 
%\end{split}
%\end{equation}
\begin{equation} \label{Firsttwocomp}
\begin{array}{lcrl}
\displaystyle{\frac{\mathrm{d}x_1}{\mathrm{d}t}} & = &A(\mathbf{x}) &\displaystyle{\frac{\partial I}{\partial x_2}} \\[4mm]
\displaystyle{\frac{\mathrm{d}x_2}{\mathrm{d}t}} & = &-A(\mathbf{x}) &\displaystyle{\frac{\partial I}{\partial x_1} }
\end{array}
\end{equation}
where $A(\mathbf{x})$ is some affine function determined by the vector field.

\begin{theorem} \label{theo-1}
The Kahan discretization of the vector field (\ref{GeneralODE}) with integral (\ref{Integral2vars}) possesses the modified integral
\begin{equation} \label{Mod2Integral}
\tilde{I}(\mathbf{x}) := \frac{I(x_1,x_2)+\frac18 h^2 D_2(a)A(\mathbf{x})^2}{1+\frac14 h^2 D_1(a)A(\mathbf{x})^2} 
\end{equation}
where
\begin{equation} \label{Determinants}
D_1(a) = \left|
\begin{array}{cc} a_1 & a_2 \\ a_2 & a_3 \end{array}
\right|,\quad
D_2(a) = \left|
\begin{array}{ccc} a_1 & a_2 & a_4\\ a_2 & a_3 & a_5\\ a_4 & a_5 & 0\end{array}
\right|.
\end{equation}
\end{theorem}

%\noindent Proof

\proof\hspace{1mm}\newline

\noindent Let the ODE
\begin{equation}\label{ODE}
\frac{d \mathbf{x}}{dt} = f(\mathbf{x}), \quad \mathbf{x} \in \mathbb{R}^n
\end{equation}
possess a quadratic integral $I$ in 2 variables. 
\begin{equation}\label{Integral}
I(x_1,x_2) = \frac{1}{2}a_1 x_1^2 + a_2 x_1x_2 +  \frac{1}{2}a_3 x_2^2 + a_4 x_1 + a_5 x_2
\end{equation}
As indicated above, it follows that the first two components of the vector field can be written
\begin{eqnarray}\label{vecfld}
\frac{dx_1}{dt} &=& A(\mathbf{x}) I_2 \\
\frac{dx_2}{dt} &=& -A(\mathbf{x}) I_1 \nonumber
\end{eqnarray}
where
\begin{equation}\label{notn1}
I_i := \frac{\partial I}{\partial x_i} \qquad i=1,2.
\end{equation}

\medskip \noindent We now discretize eqns (\ref{vecfld}) as follows:
\begin{eqnarray}\label{discvf}
\frac{x_1' - x_1}{h} &=& B(\mathbf{x} , \mathbf{x'}) I_2' + C(\mathbf{x} , \mathbf{x'}) I_2     \\
\frac{x_2' - x_2}{h} &=& -B(\mathbf{x} , \mathbf{x'}) I_1' - C(\mathbf{x} , \mathbf{x'}) I_1 \nonumber
\end{eqnarray}
where
\begin{equation}\label{notn2}
I_i' := \left. \frac{\partial I}{\partial x_i}\right|_\mathbf{x'}, \quad i=1,2.
\end{equation}
Assuming that $D_1 \neq 0$, using (\ref{Integral}), we obtain from (\ref{discvf}) that
\begin{equation}\label{equiv}
 \frac{I(x_1',x_2') - \frac{1}{2} \frac{D_2}{D_1}}{I(x_1,x_2)  - \frac{1}{2} \frac{D_2}{D_1}} 
 \equiv \frac{1+h^2 D_1 C^2(\mathbf{x} , \mathbf{x'})}{1 + h^2 D_1 B^2(\mathbf{x} , \mathbf{x'})}.
\end{equation}
Note that eq(\ref{equiv}) is an algebraic identity, obtained without any knowledge of $x_3',x_4',\dots x_n'$.

\medskip \noindent In (\ref{equiv}), the determinants $D_1$ and  $D_2$ are defined by:
\begin{eqnarray}\label{dets}
D_1 &=& a_1a_3 - a_2^2     \\
D_2 &=& 2a_2a_4a_5 - a_3a_4^2 - a_1a_5^2 \nonumber
\end{eqnarray}
We can now consider several cases:

\medskip \noindent {\bf Case(1):}
\begin{eqnarray}\label{B1C1}
B(\mathbf{x} , \mathbf{x'}) &=& E( \mathbf{x})    \\
C(\mathbf{x} , \mathbf{x'}) &=& E( \mathbf{x'})  \nonumber
\end{eqnarray}
It follows that the modified integral $\tilde{I}(\mathbf{x})$ is given by 
\begin{equation}\label{Itilde1}
\tilde{I}(\mathbf{x}) = \frac{I(x_1,x_2) - \frac{1}{2}\frac{D_2}{D_1}}{1 + h^2 D_1 E^2(\mathbf{x})}.
\end{equation}
This case includes Kahan's method for $E( \mathbf{x}) = \frac{1}{2}A(\mathbf{x})$.

\medskip \noindent {\bf Case(2):}
\begin{eqnarray}\label{B2C2}
B(\mathbf{x} , \mathbf{x'}) &=& F( \mathbf{x'})    \\
C(\mathbf{x} , \mathbf{x'}) &=& F( \mathbf{x})  \nonumber
\end{eqnarray}
It follows that the modified integral $\tilde{I}(\mathbf{x})$ is given by 
\begin{equation}\label{Itilde2}
\tilde{I}(\mathbf{x}) = \left[I(x_1,x_2) - \frac{1}{2}\frac{D_2}{D_1}\right][1 + h^2 D_1 F^2(\mathbf{x})].
\end{equation}
This case includes the trapezoidal rule for $F( \mathbf{x}) = \frac{1}{2}A(\mathbf{x})$.

\medskip \noindent {\bf Case(3):}
\begin{equation}\label{B3C3}
B(\mathbf{x} , \mathbf{x'}) = C(\mathbf{x} , \mathbf{x'}) \nonumber
\end{equation}
This case corresponds to the Discrete Gradient Method as applied to a quadratic integral. It follows that the discretization preserves the original integral $I(x_1,x_2)$.

\medskip \noindent This case has at least 2 subcases:

\medskip {\bf Case(3a):} $B(\mathbf{x} , \mathbf{x'}) = \frac{1}{2}A\left(\frac{\mathbf{x} + \mathbf{x'}}{2} \right)$ (The midpoint rule).

\medskip {\bf Case(3b):} $B(\mathbf{x} , \mathbf{x'}) = \frac{1}{2}A(\mathbf{x})$. This is an almost explicit method that we have noted before.

\bigskip \noindent {\bf Comments:}

\medskip \noindent (1) Note that the integral $\hat{I}(\mathbf{x}) = \tilde{I}(\mathbf{x}) + \frac{1}{2}\frac{D_2}{D_1}$ is given by
\begin{equation}\label{Ihat1}
\hat{I}(\mathbf{x}) = \frac{I(x_1,x_2) + \frac{1}{2}h^2D_2E^2(\mathbf{x})}{1 + h^2 D_1 E^2(\mathbf{x})} \qquad \mbox{(Case(1))} \nonumber
\end{equation}
and
\begin{equation}\label{Ihat2}
\hat{I}(\mathbf{x}) = I(x_1,x_2) + h^2 F^2(\mathbf{x}) \left( D_1 I(x_1,x_2) - \frac{1}{2} D_2  \right)       \quad \mbox{(Case(2))} \nonumber
\end{equation}
are also defined in the case $D_1=0$. We have checked that in the limit $D_1 \rightarrow 0$ these formulas are correct.

\medskip \noindent (2) We have not used anywhere that the function $A(\mathbf{x})$ should be affine. It follows that the above results remain true for any function  $A(\mathbf{x})$ (Except for the fact that in case(1) the method will not be linearly implicit if  $A(\mathbf{x})$ is not affine).

\bigskip \noindent {\bf Some Examples:}

\medskip \noindent Example 1. \cite{PPS}: 2D Suslov system
\begin{eqnarray}\label{2DSuslov}
\frac{d x_1}{dt} &=& 2\alpha x_1x_2 \\
\frac{d x_2}{dt} &=& -2x_1^2.    \nonumber
\end{eqnarray}
This system may be written
\begin{eqnarray}
\frac{d x_1}{dt} &=& 2x_1\frac{\partial I}{\partial x_2} \\
\frac{d x_2}{dt} &=& -2x_1\frac{\partial I}{\partial x_1}  \nonumber
\end{eqnarray}
with
\begin{equation}
I(x_1,x_2) = \frac{1}{2}x_1^2 + \frac{1}{2}\alpha x_2^2.     \nonumber
\end{equation}
Theorem 1 explains that the Kahan discretization of (\ref{2DSuslov}) preserves the modified integral
\begin{equation}
\tilde{I}(x_1,x_2)=  \frac{\frac{1}{2}x_1^2 + \frac{1}{2}\alpha x_2^2}{1 + h^2 \alpha x_1^2}.    \nonumber
\end{equation}

\medskip \noindent Example 2. \cite{PPS}: Zhukovsky-Volterra system with vanishing $\beta_3$
\begin{eqnarray}
\frac{d x_1}{dt} &=& \alpha x_2 x_3 -  \beta_2 x_3 \label{ZV1}   \\ 
\frac{d x_2}{dt} &=& \beta_1 x_3  \label{ZV2}    \\ 
\frac{d x_3}{dt} &=& -\alpha x_1 x_2 - \beta_1 x_2 + \beta_2 x_1  \label{ZV3} 
\end{eqnarray}
Equations (\ref{ZV1}) and  (\ref{ZV2}) may be written
\begin{eqnarray}
\frac{d x_1}{dt} &=&  x_3 \frac{\partial I}{\partial x_2}   \\ 
\frac{d x_2}{dt} &=& -x_3 \frac{\partial I}{\partial x_1}   
\end{eqnarray}
with
\begin{equation} 
I(x_1,x_2) = \frac{1}{2} \alpha x_2^2 -  \beta_1 x_1 - \beta_2 x_2.      \nonumber 
\end{equation}
Theorem 1 explains that the Kahan discretization of (\ref{ZV1}, \ref{ZV2}, \ref{ZV3}) preserves the modified integral 
\begin{equation}
\tilde{I}(x_1,x_2,x_3) = \frac{1}{2} \alpha x_2^2 - \beta_1 x_1 - \beta_2 x_2 - \frac{1}{8}h^2 \alpha \beta_1^2 x_3^2.  \nonumber
\end{equation}

\medskip \noindent Example 3. \cite{PPS}: Two coupled Euler tops
\begin{eqnarray} \label{cetvf}
\frac{dx_1}{dt} &=& \alpha_1x_2x_3   \\
\frac{dx_2}{dt} &=&  \alpha_2x_3x_1   \\
\frac{dx_3}{dt} &=&  \alpha_3x_1x_2 +  \alpha_4x_4x_5  \\
\frac{dx_4}{dt} &=&  \alpha_5x_5x_3  \\
\frac{dx_5}{dt} &=&  \alpha_6x_3x_4  
\end{eqnarray}
Equations (34) and (35) may be written
\begin{eqnarray} \label{cetvf12}
\frac{dx_1}{dt} &=& x_3 \frac{\partial I_1}{\partial x_2}  \\
\frac{dx_2}{dt} &=&  -x_3  \frac{\partial I_1}{\partial x_1}   \nonumber
\end{eqnarray}
with
\begin{equation}\label{I1}
I_1(x_1,x_2) = \frac{\alpha_1}{2}x_2^2 - \frac{\alpha_2}{2}x_1^2.
\end{equation}
Moreover (37) and (38) may be written
\begin{eqnarray} \label{cetvf45}
\frac{dx_4}{dt} &{}& = x_3 \frac{\partial I_2}{\partial x_5}  \\
\frac{dx_5}{dt} &{} & =  -x_3  \frac{\partial I_2}{\partial x_4}   \nonumber
\end{eqnarray}
with
\begin{equation}\label{I2}
I_2(x_4,x_5) = \frac{\alpha_5}{2}x_5^2 - \frac{\alpha_6}{2}x_6^2.
\end{equation}
Theorem 1 explains that the Kahan discretisation of (34-38) preserves the modified integrals
\begin{equation}\label{I1tilde}
\tilde{I}_1(x_1,x_2,x_3) =   \frac{\frac{1}{2}\alpha_1x_2^2 - \frac{1}{2}\alpha_2x_1^2}{1-\frac{h^2}{4}\alpha_1\alpha_2x_3^2} \quad ; \quad
\tilde{I}_2(x_3,x_4,x_5) =   \frac{\frac{1}{2}\alpha_5x_5^2 - \frac{1}{2}\alpha_6x_6^2}{1-\frac{h^2}{4}\alpha_5\alpha_6x_3^2}
\end{equation}
If the super-integrability condition
\begin{equation}\label{superint}
\alpha_1 \alpha_2 = \alpha_5 \alpha_6
\end{equation}
holds, eqs (34-38) may be written
\begin{eqnarray} \label{cetvfsi1}
\frac{d(x_1+x_4)}{dt} & = x_3(\alpha_1x_2 + \alpha_5x_5)  \\
\frac{d(\alpha_1x_2 + \alpha_5x_5)}{dt} & =  x_3 \alpha_1 \alpha_2 (x_1 + x_4)   \nonumber
\end{eqnarray}
Defining 
\begin{eqnarray} \label{defXY}
X &:=& x_1 + x_4  \\
Y &:=& \alpha_1x_2 + \alpha_5x_5  \nonumber
\end{eqnarray}
eq(\ref{cetvfsi1}) becomes 
\begin{eqnarray} \label{cetvfsi2}
\frac{dX}{dt} & = x_3 \frac{\partial I_3}{dY}  \\
\frac{dY}{dt} & =  -x_3 \frac{\partial I_3}{dX}   \nonumber
\end{eqnarray}
with
\begin{equation}\label{I3}
I_3(X,Y) = \frac{1}{2}Y^2 - \frac{\alpha_1 \alpha_2}{2}X^2{\color{red}.} 
\end{equation}
i.e. a quadratic function of the two variables $X$ and $Y$.

\medskip \noindent Since the Kahan discretisation is the restriction of a Runge-Kutta method to quadratic vector fields, and since all Runge-Kutta methods commute with all affine transformations (and hence with the transformation  (\ref{defXY})), Theorem 1 also explains why the Kahan discretisation preserves the modified integral
\begin{equation}\label{I3tilde}
\tilde{I}_3(x_1,x_2,x_3,x_4,x_5) =  \frac{\frac{1}{2}(\alpha_1x_2 + \alpha_5x_5)^2 - \frac{1}{2}\alpha_1\alpha_2(x_1 + x_4)^2}{1-\frac{h^2}{4} \alpha_1\alpha_2x_3^2}{\color{red}.} 
\end{equation}

\section{A family of integrable maps in $\mathbb{R}^3$}
Define
\begin{equation} \label{HandK}
\begin{array}{lcl}
H(x,y) &=&\displaystyle{ \frac12 a_1 x^2 + a_2 xy + \frac12 a_3 y^2 + a_4 x + a_5 y }\\[3mm]
K(y,z) &=& \displaystyle{\frac12 b_1 y^2 + b_2 yz + \frac12 b_3 z^2 + b_4 y + b_5 z } 
\end{array}
\end{equation}
Consider the integrable (and divergence-free) Nambu system
\begin{equation} \label{Nambu}
\morespacearray
\frac{\mathrm{d}x}{\mathrm{d}t} = \nabla H \times \nabla K = \left(
\begin{array}{r}
\frac{\partial K}{\partial z} \frac{\partial H}{\partial y}  \\
-\frac{\partial K}{\partial z} \frac{\partial H}{\partial x}  \\
\frac{\partial H}{\partial x} \frac{\partial K}{\partial y} 
\end{array}
\right)
\end{equation}

\begin{theorem}\label{theo-2}
The Kahan discretization of vector field (\ref{Nambu}) with integrals (\ref{HandK}) possesses the modified integrals
\begin{eqnarray}
\tilde{H}(x,y,z) &=& \frac{H(x,y)+\frac18 h^2 D_2(a)(b_2y+b_3z+b_5)^2}{1+\frac14h^2 D_1(a)(b_2y+b_3z+b_5)^2} \label{Ht} \\
\tilde{K}(x,y,z) &=& \frac{K(y,z)+\frac18 h^2 D_2(b)(a_1x+a_2y+a_4)^2}{1+\frac14h^2 D_1(b)(a_1 x+a_2y+a_4)^2} \label{Kt}
\end{eqnarray}
and preserves the modified measure
\begin{equation} \label{Gmeas}
   g(x,y,z)\,  \mathrm{d}x\wedge \mathrm{d}y\wedge \mathrm{d}z
\end{equation}
with
\begin{equation} \label{Gdef}
g(x,y,z) = \big(1+\frac14 h^2 D_1(a)(b_2y + b_3 z + b_5)^2\big)^{-1}
\big(1+\frac14 h^2 D_1(b)(a_1x + a_2 y + a_4)^2\big)^{-1}
\end{equation}
where $D_1(a)$, $D_2(a)$ are defined in (\ref{Determinants}) and $D_1(b)$, $D_2(b)$ similarly. It follows that the Kahan discretization of (\ref{Nambu}) is completely integrable.
\end{theorem}
\proof\hspace{1mm}\newline
\renewcommand{\labelenumi}{(\roman{enumi})}
\begin{enumerate}
\item  The preservation of the two integrals (\ref{Ht}) and (\ref{Kt}) follows using Theorem~\ref{theo-1} and
\begin{equation} \label{Kz-Hx}
\frac{\partial K}{\partial z} = b_2 y + b_3 z + b_5,\quad\mbox{resp}\ \frac{\partial H}{\partial x} = a_1 x + a_2 y + a_4
\end{equation}
\item By definition, the measure $g\,\mathrm{d}x\wedge \mathrm{d}y\wedge \mathrm{d}z$ is preserved if
\begin{equation}\label{intcond}
\int g(x,y,z)\, \mathrm{d}x\wedge \mathrm{d}y\wedge \mathrm{d}z
=\int g(x',y',z')\,\mathrm{d}x'\wedge \mathrm{d}y'\wedge \mathrm{d}z',
\end{equation}
hence if the so-called density $g$ satisfies
\begin{equation} \label{DetCondition}
\morespacearray
g(x,y,z) = g(x',y',z')
\left|
\begin{array}{ccc}
\frac{\partial x'}{\partial x} &  \frac{\partial x'}{\partial y} &  \frac{\partial x'}{\partial z}\\
\frac{\partial y'}{\partial x} &  \frac{\partial y'}{\partial y} &  \frac{\partial y'}{\partial z} \\
\frac{\partial z'}{\partial x} &  \frac{\partial z'}{\partial y} &  \frac{\partial z'}{\partial z} 
\end{array}
\right|.
\end{equation}
In the case at hand, the map $(x,y,z)\mapsto (x',y',z')$ is given by Kahan's map, and the identity (\ref{DetCondition}) has been verified using Maple, after substituting (\ref{Gdef}).
\item By definition, integrability of a three-dimensional map follows directly from the preservation of two integrals plus a preserved measure.
\end{enumerate}

\section*{Comments:}
\begin{enumerate}
\item We note that the functional form of the preserved integrals (\ref{Ht}) and (\ref{Kt}), and density (\ref{Gdef}) is not unique, because any function of the integrals is an integral, and the product of the density with any integral will be a preserved density.

\medskip In particular, if $D_1(a)$ and $D_1(b)$ do not vanish, alternative discrete integrals  are given by
\begin{eqnarray}
\hat{H}(x,y,z) &=& \frac{H(x,y) - \frac{1}{2}\frac{D_2(a)}{D_1(a)}}{1+\frac14h^2 D_1(a)(b_2y+b_3z+b_5)^2}  \\
\hat{K}(x,y,z) &=& \frac{K(y,z) - \frac{1}{2}\frac{D_2(b)}{D_1(b)}}{1+\frac14h^2 D_1(b)(a_1 x+a_2y+a_4)^2},
\end{eqnarray}
and an alternative preserved density is given by
\begin{equation}\label{alt_g}
g(x,y,z) = \left( H(x,y) - \frac{1}{2}\frac{D_2(a)}{D_1(a)} \right)^{-1} \left( K(y,z) - \frac{1}{2}\frac{D_2(b)}{D_1(b)} \right)^{-1}.
\end{equation}

\item Note that the density (\ref{alt_g}) does not depend on the timestep $h$, and therefore is also preserved by the ODE (\ref{Nambu}).

\item Of course, if the reader so chooses, it is possible to introduce normal forms for this family of maps by applying appropriate affine transformations to the coordinates.

\item The ODE (\ref{Firsttwocomp}) is invariant under $I \rightarrow \alpha I$, $A \rightarrow \beta A$, $t \rightarrow \frac{t}{\alpha \beta}$. Similarly, the modified integral (\ref{Mod2Integral}) preserved by the Kahan discretization is covariant under $I \rightarrow \alpha I$, $A \rightarrow \beta A$, $h \rightarrow \frac{h}{\alpha \beta} I$.  

\item For some examples of other integrable families of maps in $\mathbb{R}^3$ published in the literature, the reader is referred to \cite{CMMOQ,I,RIQ,HKY,MT,FK}. The maps in \cite{CMMOQ} are closest to the maps in the current paper. Nevertheless they are different: the integrals of the ODE in \cite{CMMOQ} are essentially homogeneous, whereas the integrals of the ODE (\ref{Nambu}) are generically inhomogeneous.

\end{enumerate}

\section*{Acknowledgements}
This research was supported by the Australian Research Council and by the European Unions Horizon 2020
research and innovation programme under the Marie Sk\l{}odowska-Curie grant
agreement No. 691070, and by The Research Council of Norway. GRWQ is grateful to Jason Frank for valuable discussions, and to Khaled Hamad for Maple assistance with the proof of Theorem~\ref{theo-1}. 

%The computations in this paper were performed by using Maple\texttrademark.

\end{document}